\documentclass[11pt,a4paper]{article}
 
\usepackage[english]{babel}
\usepackage{amsmath,amsthm,amssymb,graphicx,xcolor,bm}  
\usepackage{multirow}

 \allowdisplaybreaks

\newtheorem{theorem}{Theorem}
\newtheorem{lemma}{Lemma}

\theoremstyle{definition}
\newtheorem{remark}{Remark}
 \newtheorem{example}{Example}

\begin{document}

 \title{Pinsker bound under measurement budget constrain: optimal allocation} 
\author{Eduard Belitser \\ {\it VU Amsterdam}}
\date{}   
\maketitle

\renewcommand{\abstractname}{\vspace{-2\baselineskip}}
\begin{abstract}
In the classical many normal means with different variances, we consider
the situation when the observer is allowed to allocate the available 
measurement budget over the coordinates of the parameter of interest. 
The benchmark is the minimax linear risk over a set.
We solve the problem of optimal allocation of observations under the measurement 
budget constrain for two types of sets, ellipsoids and hyperrectangles.
By elaborating on the two examples of Sobolev ellipsoids and hyperectangles, we  
demonstrate how re-allocating the measurements in the (sub-)optimal way 
improves on the standard uniform allocation. In particular, 
we improve the famous Pinsker (1980) 
bound. 
\end{abstract}

\vspace{-2\baselineskip} 
\footnote{\vspace{-\baselineskip}\\
\textit{MSC2010 subject classification:}
primary 62C20, 62G05; secondary 62G20. }

\section{Introduction}
Suppose we observe independent
$X_{ij} \sim N(\theta_i, \sigma_i^2)$, $j=1,\ldots n_i$, $i \in \mathbb{N}$, 
where $\sigma_i\ge \sigma>0$ (without loss of generality assume $\sigma=1$) 
and the goal is to recover the unknown parameter 
$\theta=(\theta_i,\, i \in \mathbb{N})\in \Theta$, for some 
$\Theta \subseteq \ell_2$ to be introduced later. 
By sufficiency, this setting is equivalent to the the classical many 
normal means model: with $X_i = n_i^{-1}\sum_{j=1}^{n_i} X_{ij}$, 
\begin{align}
\label{model}  
X_i \sim N(\theta_i, \sigma_i^2/n_i), \quad  i \in \mathbb{N}.
\end{align}
From now on we study model \eqref{model} and allow 
$n_i \in \mathbb{R}_+=\{a\in\mathbb{R}: \, a \ge 0\}$.

Let us call the vector $\bar{n} = (n_i,\, i \in \mathbb{N})$ 
\emph{measurement allocation} 
and $B(\bar{n}) = \sum_i n_i$ its \emph{measurement budget}, which may be infinite.
Introduce the (quadratic) minimax risk 
$r(\bar{n},\Theta)=\inf_{\hat{\theta}}\sup_{\theta \in \Theta} 
\mathrm{E}_\theta \|\hat{\theta}-\theta\|^2$, where the infimum is taken 
over all possible estimators $\hat{\theta}=\hat{\theta}(X)$, measurable 
functions of the data $X=(X_i, i \in \mathbb{N})$. 
Here $\|\theta\| = (\sum_i \theta_i^2)^{1/2}$ is the $\ell_2$-norm and  
$\sum_i$ means the summations over $\mathbb{N}$ from now on. 
The classical situation is the so called \emph{uniform measurement allocation}
$\bar{n}_u(n)=\bar{1}n=
(n,n,\ldots)$, where $\bar{1}=(1,1,\ldots)$. 
This is a well studied setting, typically under the asymptotic regime as $n \to \infty$.  
In his seminal work, Pinsker (1980) studied this case  and derived the exact asymptotic 
behavior of the minimax risk $r(\bar{n}_u(n),\Theta)$ 
for ellipsoidal sets $\Theta = \mathcal{E}$; see the exact definition below.
For $\sigma_i=\sigma$, many extensions and various aspects of the uniform 
allocation framework have been investigated by Donoho et al.\ (1990). 
If $\sigma_i \to \infty$ as $i\to \infty$ in \eqref{model}, the problem of  recovering $\theta$  
is known to be the so-called \emph{inverse problem}.

In this note, the measurement allocation $\bar{n}$ is not necessarily uniform:
$n_i$ may vary with  $i$. In fact, the observer is allowed to allocate the measurements 
in any way, but under the so called \emph{measurement budget constrain}:
\begin{align}
\label{mbc}
\bar{n}  
\in\mathcal{N}_b = \big\{\bar{m}=(m_i,\, i \in \mathbb{N}): \, 
m_i \ge 0,\, B(\bar{m})=\sum\nolimits_i m_i\le b\big\},
\end{align}
where $n_i$ has the interpretation of number of measurements of $\theta_i$.

The goal of this note  is twofold: first,  to derive the optimal (in a certain sense explained 
in the next section) allocation $\bar{n}_o\in\mathcal{}N_b$ under the measurement 
budget constrain for two families of sets $\Theta\subseteq \ell_2$, hyperrectangles 
and ellipsoids; second, to demonstrate that re-allocating measurements 
in the (sub-)optimal way improves upon the standard uniform allocation scheme, 
and to quantify this improvement. Some inequalities for the beta function
are obtained as consequence.

The paper is organized as follows. 
In the next section we present the exact mathematical framework 
for addressing our goals. In Section \ref{sec_ellipsoids}, we derive an implicit description  
of the optimal allocation for ellipsoids. We also find an explicit sub-optimal solution 
and elaborate on this for the example of Sobolev ellipsoids. 
In particular, we demonstrate the  improvement  
of the famous Pinsker bound provided by re-allocating the measurements 
in the proposed sub-optimal way. This gives a bound on the performance of
the optimal (implicit) allocation. In Section \ref{sec_hyperrectangles}, 
we solve the problem of optimal allocation for general hyperrectangles 
and present the exact evaluation for Sobolev hyperrectangles. 
Some peculiar inequalities for the beta function are obtained as consequence.
A short discussion of possible extensions is given in Section \ref{sec_discussion}.

\section{Preliminaries and statement of the problem} 
\label{sec_prelim}

Introduce the class of the linear estimator 
$\hat{\theta}(\lambda) =(\lambda_i X_i, \, i \in \mathbb{N})$
with the weights $\lambda=( \lambda_i,\, i \in \mathbb{N})$.
Our benchmark is the minimax linear risk: 
\begin{align}
\label{lin_minimax}
\mathnormal{R}(\bar{n},\Theta)=
\inf_\lambda \sup_{\theta \in \Theta} \mathrm{E}_\theta \|\hat{\theta}(\lambda) - \theta\|^2=
\inf_\lambda\sup_{\theta \in \Theta} \sum_i \Big(\frac{\sigma^2_i\lambda_i^2}{n_i}
+(1-\lambda_i)^2 \theta_i^2 \Big).
\end{align} 
In this note, we focus on the minimax linear risk instead of 
the minimax risk $r_n(\bar{n},\Theta)=\inf_{\hat{\theta}}\sup_{\theta \in \Theta} 
\mathrm{E}_\theta \|\hat{\theta}-\theta\|^2$ for the following reasons. 
Firstly, as we will see below, for some important sets 
$\Theta$ (here we consider hyperrectangles and 
ellipsoids) the minimax linear risk 
$\mathnormal{R}(\bar{n},\Theta)$ is a more tractable quantity than 
the minimax risk $r(\bar{n},\Theta)$. Secondly, for the considered 
choices of set $\Theta$ (hyperrectangles and ellipsoids) we have the bound
\begin{align}
\label{const_C_Theta}
 r(\bar{n},\Theta)\le R(\bar{n},\Theta) 
\le C_\Theta r(\bar{n},\Theta),
\end{align}
for some absolute constant $C_\Theta \ge 1$, so that the minimax rate,
in the asymptotic regime $n \to \infty$, is the same for the both risks.

As for set $\Theta$,  we introduce general ellipsoids and hyperrectangles 
\begin{align}
\label{ellipsoids}
\mathcal{E}(a) = \big\{ \theta \in \ell_2: \, \sum\nolimits_i (\tfrac{\theta_i}{a_i})^2 \le 1\big\}, \quad
\mathcal{H}(a)= \{ \theta \in \ell_2: \, |\theta_i|\le a_i, \, i \in \mathbb{N}\},
\end{align}
where $a=(a_i,\,  i \in \mathbb{N})$ is a sequence of numbers in 
$[0,+\infty]$ which converge to 0 as $i\to \infty$.
We adopt the conventions $0/0=0$ and $x/(+\infty)=0$ for $x \in \mathbb{R}$. 
Without loss of generality, let the sequence $a_i$ be strictly positive 
and nonincreasing. 

In case of the uniform allocation and constant $\sigma_i=\sigma$,
Donoho et al.\ (1990) established that $C_\Theta=1.25$ (can be replaced by $1.247$) 
in \eqref{const_C_Theta} for a compact, orthosymmetric, convex and  quadratically convex set 
$\Theta \subset \ell_2$. If we now turn to the case of arbitrary allocation $\bar{n}$ and 
arbitrary positive $\sigma_i$'s, the same bound immediately follows from  Donoho et al.\ (1990) 
for hyperrectangles. Moreover, as is shown by Reshetov (2010),
\eqref{const_C_Theta} still holds with $C_\Theta=1.25$  for any compact 
(can be extended to certain noncompact cases), 
orthosymmetric, convex and quadratically convex sets of $\ell_2$.
Such sets include hyperrectangles, ellipsoids and $\ell_p$-bodies with $p\ge 2$.
For the standard uniform allocation $\bar{n}=\bar{n}_u(n)$, 
under some mild conditions (excluding pathological ellipsoids  $\mathcal{E}$) 
Pinsker (1980) established  a remarkable result $r(\bar{n}_u(n),\mathcal{E}) 
=R(\bar{n}_u(n),\mathcal{E}) (1+o(1))$ as $n \to \infty$, cf.\ Nussbaum (1996).

The \emph{optimal allocation problem over set $\Theta$}:
with $\mathcal{N}_n$ defined by \eqref{mbc}, 
\begin{align}
\label{opt_allocation}
R(\bar{n}_o, \Theta) = \inf_{\bar{n} \in \mathcal{N}_n} R(\bar{n},\Theta).
\end{align}
The goal is to derive the optimal measurement allocation 
$\bar{n}_o=\bar{n}_o(n)=\bar{n}_o(n,\Theta)$  according to \eqref{opt_allocation} and 
to study the asymptotic behavior of $R(\bar{n}_o, \Theta)$ as $n \to \infty$. 
Henceforth, the default asymptotic regime (e.g., for $o(1)$) is $n\to \infty$,
unless otherwise specified.
For positive $g_n,h_n$, 
$g_n\sim h_n$ means $g_n=h_n(1+o(1))$, $g_n \asymp h_n$ 
means that $g_n=O(h_n)$ and $h_n =O(g_n)$.

In this note, we solve the problem \eqref{opt_allocation} by determining the optimal allocation 
$\bar{n}_o(n,\Theta)$ and the corresponding optimized risk $R(\bar{n}_o, \Theta)$
for ellipsoids $\Theta=\mathcal{E}(a))$ and hyperrectangles 
$\Theta=\mathcal{H}(a))$. For the two particular 
examples, Sobolev ellipsoids and  Sobolev hyperrectangles 
in the mildly ill-posed model, we derive sharp asymptotic expressions.

We also show that by re-allocating measurements one can improve on 
the standard uniform allocation. A way to quantify the amount of improvement
of the uniform allocation $\bar{n}_u$ by a contender (re-)allocation 
$\bar{m}= (m_i,\, i \in \mathbb{N})$ (e.g., $\bar{m}=\bar{n}_o$, 
the optimal allocation \eqref{opt_allocation})
is by relating their risks for $\bar{n}_u$.
However, relating, for example, the risks for $\bar{n}_u(n)$ and $\bar{n}_o(n)$ 
is not a fair comparison between allocations $\bar{n}_u(n)$ and $\bar{n}_o(n)$.
The risk for any parsimonious allocation $\bar{n}\in\mathcal{N}_n$ 
(including $\bar{n}_o(n)$) 
will expectedly  be bigger than the risk for the generous uniform allocation 
scheme $\bar{n}_u(n)$ which allows $n$ measurements for each coordinate 
of infinite dimensional $\theta$. We will see this for Sobolev hyperrectangles and ellipsoids.  

To compare in a fair way different allocations, their measurement budgets 
must be matched. 
However, an allocation $\bar{n}$ can have an infinite measurement budget; 
e.g., $\bar{n}_u(m)$, and in general $\bar{n}_u(m)$ is not matchable 
(in terms of measurement budgets) to any $\bar{n} \in \mathcal{N}_n$ 
for any $m,n\ge 0$.
The following framework tackles this.

For allocation $\bar{n}$ and a \emph{pattern} $\bar\delta=(\delta_i,\, i \in \mathbb{N})$ with 
$\delta_i\in \{0,1\}$, define the $\bar\delta$-pattern of $\bar{n}$ as the entrywise 
product $\bar{n} \cdot \bar\delta=(\delta_i n_i, \, i\in\mathbb{N})$. 
Introduce the set of \emph{effective patterns}
$
\mathcal{D}_e(\bar{n})=\mathcal{D}_e(\bar{n},\Theta)=
\{\delta: R(\bar{n},\Theta) = R(\bar{n}\cdot \bar{\delta}, \Theta)\}.
$
This set is not empty as $\bar{1} \in \mathcal{D}_e(\bar{n})$.
To assess the improvement of  an arbitrary allocation $\bar{n}$ by using  
a contender allocation $\bar{m}$  with 
$B(\bar{m})<\infty$, consider two cases. 

Case I: for some pattern $\bar\delta\in\mathcal{D}_e(\bar{n})$,
$B(\bar{n}\cdot \bar{\delta})=\sum_i n_i \delta_i < \infty$.  In this case, 
the so called \emph{effective measurement budget} (EMB) 
$B_{e}(\bar{n})=B_{e}(\bar{n},\Theta)=
\inf_{\delta\in\mathcal{D}_e(\bar{n})} B(\bar{n} \cdot \bar{\delta})$ 
of $\bar{n}$ is finite, so that the EMB's  
of $\bar{n}$ and $\bar{m}$ can be matched: $B(\bar{m})=B_e(\bar{n})$. 
The amount of improvement of allocation $\bar{n}$ 
by reallocation $\bar{m}$ is then measured by the risk ratio
\begin{align}
\label{risk_ratio}
\rho_1(\bar{n},\bar{m},\Theta)=\frac{R(\bar{n},\Theta)}{R(\bar{m},\Theta)},
\quad \text{with} \;\; B(\bar{m})=B_e(\bar{n},\Theta).
\end{align}

Case II: there is no $\bar\delta\in\mathcal{D}_e(\bar{n})$ such that
$B(\bar{n}\cdot \bar{\delta})=\sum_i n_i \delta_i < \infty$.  In this case, 
we quantify the amount of improvement of $\bar{n}$ by $\bar{m}$ with 
$B(\bar{m})<\infty$ as follows:
with $\mathcal{B}(\bar{n},\bar{m})  
=\{\bar{\delta}: B(\bar{n}\cdot\bar{\delta})\le B(\bar{m})\}$,
\begin{align}
\label{risk_ratio1}
\rho_2(\bar{n},\bar{m},\Theta)=\frac{\inf_{\bar{\delta}\in \mathcal{B}(\bar{n},\bar{m})} 
R(\bar{n}\cdot \bar{\delta}, \Theta)}{R(\bar{m},\Theta)}.
\end{align}
The numerator of the ratio 
is the best risk over all $\bar{\delta}$-patterns 
of $\bar{n}$ with the bounded (by $B(\bar{m})$) budget. We consider  only
uniform allocations $\bar{n}=\bar{n}_u(k)$. In this case, the numerator in \eqref{risk_ratio1}
will be optimized also with respect to the parameter $k$: $\inf_{\bar{\delta}, k: \, 
\bar{\delta}\in \mathcal{B}(\bar{n}_u(k),\bar{m})}R(\bar{n}_u(k)\cdot \bar{\delta}, \Theta)$.
 
Typically, one would like to compare with the optimal reallocation 
\eqref{opt_allocation}: $\bar{m} = \bar{n}_o(n)$. 
Notice that in this case we can reduce $\rho_1(\bar{n},\bar{n}_o,\Theta)$ to 
$\rho_2(\bar{n},\bar{n}_o,\Theta)$ by re-parametrizing the budgets.  
Clearly, always $\rho_i(\bar{n},\bar{n}_o,\Theta) \ge 1$, $i=1,2$, and 
the bigger $\rho_i$, the bigger the improvement.

\section{Optimal allocation for ellipsoids}
\label{sec_ellipsoids}

Consider the ellipsoidal set $\Theta=\mathcal{E}=\mathcal{E}(a)$ defined in \eqref{ellipsoids}. 
Introduce the risk of the linear estimator $\hat{\theta}(\lambda)$ 
at point $\theta\in \ell_2$: 
$R(\bar{n},\lambda,\theta) = \mathrm{E}_\theta 
\|\hat{\theta}(\lambda) - \theta\|^2=\sum_i \big[\frac{\sigma^2_i\lambda_i^2}{n_i}
+(1-\lambda_i)^2 \theta_i^2 \big]$.
Recall the minimax linear risk \eqref{lin_minimax} over $\mathcal{E}$:
$\mathnormal{R}(\bar{n},\mathcal{E})=\inf_\lambda\sup_{\theta \in \mathcal{E}} 
R(\bar{n},\lambda,\theta)$.
The following technical lemma from Belitser and Levit (1995) 
describes  $\mathnormal{R}(\bar{n},\mathcal{E})$. 
For $b \in \mathbb{R}$, denote $b_+=\max\{b,0\}$.

\begin{lemma} 
\label{lem1}
\begin{align}
\mathnormal{R}(\bar{n},\mathcal{E})
&=\sup_{\theta \in \mathcal{E}} \inf_\lambda
R(\bar{n},\lambda,\theta) 
=R(\bar{n},\lambda_o,\theta_o)\notag\\
&=\sup_{\theta \in \mathcal{E}}\sum_i\frac{\theta_i^2\sigma_i^2}
{n_i \theta_i^2+\sigma_i^2}=\sum_i \frac{\sigma_i^2}{n_i} (1- t a_i^{-1})_+,
\label{d_formula}
\end{align}
where  the saddle point $\lambda_{o,i}= (1-ta_i^{-1})_+$,  
$\theta^2_{o,i}=\frac{\sigma^2_i
a_i (1- ta_i^{-1})_+}{n_i t}$ and $t=t(\bar{n})=t(\bar{n},\mathcal{E})$ 
is the solution of the equation (with the conventions $0/0=0$)
\begin{align}
\label{t_E}
\sum\nolimits_i \frac{\sigma^2_i(1-ta_i^{-1})_+}{n_i a_i} =t. 
\end{align}
\end{lemma}
 
Lemma \ref{lem1} also determines the effective measurement budget of $\bar{n}$:
\begin{align}
\label{EMB}
B_{e}(\bar{n}) =\sum_{i=1}^{d(\bar{n})} n_i, \;\; \text{where}\;\;
d(\bar{n})=d(\bar{n},\mathcal{E})= \max\{k\in\mathbb{N}: \, a_k < t(\bar{n})\}
\end{align}
has the meaning of the \emph{number of effectively 
estimated coordinates} of $\theta\in\mathcal{E}$.
Indeed, for each $\theta \in \Theta$ the optimal (in terms of the risk $R$) 
estimator does not spend any budget on measuring the coordinates 
$\theta_k$, $k>d(\bar{n})$. 

In view of \eqref{d_formula}, 
the optimal allocation problem \eqref{opt_allocation} is formally solved 
by the following theorem.
\begin{theorem}
\label{th1}
Let $t=t(\bar{n})$ be defined by \eqref{t_E}.
The optimal allocation $\bar{n}_o=\bar{n}_o(n, \mathcal{E})$
and the optimized risk  $R(\bar{n}_o,\mathcal{E})$ are
determined by 
\[
R(\bar{n}_o,\mathcal{E}) =  
\inf_{\bar{n} \in \mathcal{N}_n} \sum_i \frac{\sigma_i^2}{n_i} 
(1- t(\bar{n}) a_i^{-1})_+.
\]
\end{theorem}
Unfortunately, this theorem does not provide explicit analytic formulas for 
the optimal allocation $\bar{n}_o(n,\mathcal{E})$ and the optimized risk  
$R(\bar{n}_o,\mathcal{E})$. Below we provide some \emph{sub-optimal
solution} that has a more explicit form and is still intended to improve on 
the uniform allocation. This sub-optimal solution also gives a lower 
bound on the performance of the optimal (implicit) allocation $\bar{n}_o$.

Treating $t$ as fixed and minimizing the last expression for 
$R(\bar{n},\mathcal{E})$ in \eqref{d_formula} with respect to 
$(n_i,i\in \mathbb{N}) \in\mathcal{N}_n$ (i.e., under the restriction $\sum_i n_i =n$)
by Lagrange multiplier yields $\bar{n}_s=\bar{n}_s(n)=\bar{n}_s(n,t)$:
\begin{align}
\label{n_s}
n_{s,i}(n,t)= \frac{n\sigma_i \big[1-t a_i^{-1}\big]_+^{1/2}}
{\sum_k \sigma_k \big[1-ta_k^{-1}\big]_+^{1/2}}, 
\quad  i \in \mathbb{N}.
\end{align}
According to Lemma \ref{lem1}, $t$ must be the solution of 
\eqref{t_E} with $\bar{n}_= \bar{n}_s$ in order for the risk 
$R(\bar{n}_s,\mathcal{E})$ to be equal to the expression \eqref{d_formula}. 
Plugging in the expression \eqref{n_s} into \eqref{t_E} gives the equation 
for $t_s=t_s(n)$:
\begin{align}
\label{t_E(n)}
\sum\nolimits_k \sigma_k \big[1-t_sa_k^{-1}\big]_+ ^{1/2}  
\sum\nolimits_i \sigma_i a_i^{-1}  \big[1-t_sa_i^{-1}\big]_+ ^{1/2}=n t_s.
\end{align}
Thus the sub-optimal allocation 
$\bar{n}_s(n,t_s)$ is 
defined by \eqref{n_s} and \eqref{t_E(n)}.
The corresponding sub-optimal risk is obtained by
substituting $\bar{n}_s(n)=\bar{n}_s(n,t_s)$ into the expression 
\eqref{d_formula} for $R(\bar{n},\mathcal{E})$:
\begin{align}
\label{R_ellips}
\mathnormal{R}(\bar{n}_s(n),\mathcal{E})=
n^{-1}\Big(\sum\nolimits_i\sigma_i \big[1-t_s a_i^{-1}\big]_+^{1/2}
\Big)^2,
\end{align}
with $t_s$ defined by \eqref{t_E(n)}.
To summarize, the fomulas \eqref{n_s}--\eqref{R_ellips}
describe a sub-optimal solution to the 
measurement allocation problem over ellipsoids 
under measurement budget constrain. 
In the next example we compute the asymptotic behavior of 
this sub-optimal solution in case of Sobolev ellipsoids.

\begin{example} 
\label{example1}
Consider the Sobolev ellipsoid 
$\mathcal{E}=\mathcal{E}(a)$ with $a_i^2=Qi^{-2\alpha}$, 
$\sigma_i^2=\sigma^2 i^{2\beta}$,  $\alpha,Q,\sigma^2>0$, $\beta>-\frac{1}{2}$. 
Without loss of generality, let $\sigma^2=Q=1$.
Indeed, denote $R_\mathcal{E}(\bar{n}, C, c) = R(\bar{n}, \mathcal{E})$ 
for $a^2_i=C\tilde{a}_i^2$, $\sigma^2_i=c\tilde{\sigma}_i^2$, 
with some fixed $\tilde{a}_i$,  $\tilde{\sigma}^2_i$.
Then by Lemma \ref{lem1}, one can show that 
$R_\mathcal{E}(\bar{n},C,c) = CR_\mathcal{E}(\bar{n}C/c,1,1)$.

The following asymptotic identity holds for any $\alpha>0$, $\beta, \kappa>-1$:
\[
\sum\nolimits_k k^\beta \big(1-\tfrac{k^\alpha}{M}\big)_+^\kappa \sim
 \alpha^{-1}B\big(\tfrac{\beta+1}{\alpha}, \kappa+1\big) M^{\frac{\beta+1}{\alpha}}, 
 \quad M \to \infty,
\]
where $B(\cdot, \cdot)$ is the beta function.
Using this relation in some tedious computations gives the (asymptotic) 
solution of \eqref{t_E(n)} for the Sobolev ellipsoid:
\[
t_s=t_s(n)\sim d_s^{-\alpha}, \quad
d_s \sim \bigg(
\frac{\alpha^2(3\alpha+2\beta+2)}
{2(\beta+1)B^2(\tfrac{\beta+1}{\alpha},\tfrac{3}{2})} n\bigg)^{\frac{1}{2(\alpha+\beta+1)}},
\]
where the quantity
$d_s=d_s(n)=d_s(n,\mathcal{E})=d(\bar{n}_s(n),\mathcal{E})=
\max\{k\in\mathbb{N}: \, a_k < t_s\}$, with $t_s$ defined by \eqref{t_E(n)},
is the number of nonzero coordinates in the sub-optimal allocation vector $\bar{n}_s(n)$.
The expressions for the sub-optimal allocation and the corresponding 
sub-optimal risk follow from \eqref{n_s} and \eqref{R_ellips}:
\begin{align}
n_{s,i}(n)&=\frac{n\sigma_i \big[1-t_s a_i^{-1}\big]_+^{1/2}}
{\sum_k \sigma_k \big[1-t_sa_k^{-1}\big]_+^{1/2}}
\sim \frac{\alpha( \tfrac{i}{d_s})^\beta\big[1- (\tfrac{i}{d_s})^\alpha\big]_+^{1/2}}
{B(\tfrac{\beta+1}{\alpha},\tfrac{3}{2})} \frac{n}{d_s},
\quad  i \in \mathbb{N}, \notag\\
\label{risk_opt_alloc}
R(\bar{n}_s(n),\mathcal{E}) &=
n^{-1}\Big(\sum\nolimits_i\sigma_i \big[1-t_ s a_i^{-1}\big]_+^{1/2}\Big)^2
\sim B_{\mathcal{E}} n^{-\frac{\alpha}{\alpha+\beta+1}},
\end{align}
with $B_{\mathcal{E}} =
B_{\mathcal{E}}(\alpha,\beta)
= \big(B^2(\frac{\beta+1}{\alpha},\frac{3}{2})/\alpha^2\big)^{\frac{\alpha}{\alpha+\beta+1}}
\big(\frac{3\alpha+2\beta+2}{2(\beta+1)}\big)^{\frac{\beta+1}{\alpha+\beta+1}}$.

Using the relation 
$\sum_{k=1}^M k^\kappa \sim \frac{M^{k+1}}{\kappa+1}$
as $M \to \infty$ for $\kappa >-1$, we compute 
the risk for the uniform allocation $\bar{n}_u(n)$
(cf.\ Belitser and Levit (1995)):
\begin{align}
\label{risk_unif_alloc}
R(\bar{n}_u(n),\mathcal{E}) =n^{-1} \sum\nolimits_i 
\sigma_i^2(1-t(\bar{n}_u(n))a_i^{-1})_+ \sim
\bar{B}_{\mathcal{E}} n^{-\frac{2\alpha}{2\alpha+2\beta+1}},
\end{align}
where
$
\bar{B}_{\mathcal{E}}=\bar{B}_{\mathcal{E}}(\alpha,\beta)=
\frac{(2\alpha+2\beta+1)^{\frac{2\beta+1}{2\alpha+2\beta+1}}}
{2\beta+1} 
\big(\frac{\alpha}{\alpha+2\beta+1}
\big)^{\frac{2\alpha}{2\alpha+2\beta+1}}.
$
For $\beta=0$, $\bar{B}_{\mathcal{E}}(\alpha,0)$ is known to be the famous
\emph{Pinsker constant}; cf.\ Pinsker (1980) and Nussbaum (1996).

Note that \eqref{risk_opt_alloc} and \eqref{risk_unif_alloc} imply that
\[
R(\bar{n}_s(n),\mathcal{E})\asymp n^{-\frac{\alpha}{\alpha+\beta+1}} 
\gg n^{-\frac{2\alpha}{2\alpha+2\beta+1}} \asymp R(\bar{n}_u(n),\mathcal{E}).
\] 
This is of course expected as the sub-optimal allocation   
under measurement budget constrain $\bar{n}_s(n)$ is a 
parsimonious regime (only $n$ measurements for 
all coordinates of $\theta$ are allowed) whereas the 
uniform allocation scheme allows $n$ measurements for each 
coordinate of (infinite dimensional) $\theta$.

To illustrate how reallocating the measurements in the sub-optimal way 
\eqref{n_s} improves the uniform allocation scheme and to quantify 
this improvement, we use the risk ratio \eqref{risk_ratio} with $\bar{n}=\bar{n}_u(n)$
and $\bar{m}=\bar{n}_s(B_e(\bar{n}_u(n)))$; i.e.,
the effective measurement budgets of the both allocations are matched. 
According to \eqref{EMB}, $B_e(\bar{n}_u(n))=d_n n$, where 
$d_n=d(\bar{n}_u(n))$ is as follows (cf.\ Belitser and Levit (1995)):
\begin{align}
\label{d_n}
d_n=d_n(\mathcal{E})\sim D_\mathcal{E} n^{\frac{1}{2\alpha+2\beta+1}},\quad
{\textstyle D_\mathcal{E}=\big(\frac{(2\alpha+2\beta+1)
(\alpha+2\beta+1)}{\alpha}\big)^{\frac{1}{2\alpha+2\beta+1}}}.
\end{align}
Using \eqref{risk_opt_alloc}, \eqref{risk_unif_alloc}, \eqref{d_n} and 
performing some tedious computations, we derive the limit risk ratio 
\begin{align*}
&\rho_1(\bar{n}_u,\bar{n}_s,\mathcal{E})= 
\frac{R(\bar{n}_u(n),\mathcal{E})}{R(\bar{n}_s(B_e(\bar{n}_u(n))),\mathcal{E})}
\sim\frac{\bar{B}_\mathcal{E} D_{\mathcal{E}}^{\frac{\alpha}{\alpha+\beta+1}}}
{B_\mathcal{E}} 
=\rho_{\mathcal{E}}(\alpha,\beta)
\\
&= \frac{1}{2\beta+1} 
\Big(\frac{2(\beta+1)(2\alpha+2\beta+1)}{3\alpha+2\beta+2}\ \Big)^{\frac{\beta+1}{\alpha+\beta+1}}
 \Big(\frac{\alpha^3B^{-2}(\frac{\beta+1}{\alpha}, \frac{3}{2})}{\alpha+2\beta+1}\Big)^{\frac{\alpha}{\alpha+\beta+1}}.
\end{align*}
Table \ref{table:1} presents a small selection of  the computed limit risk ratio  
$\rho_{\mathcal{E}}(\alpha,\beta)$
for several values of parameters $\alpha$ and $\beta$.
\begin{table}[ht]
\centering
\begin{tabular}{| c | c  c c c c c c c |} 
\hline
\multirow{2}*{$\beta$} 
& \multicolumn{8}{ c|}{$\alpha$}   \\ \cline{2-9}                
 & $0.5$ & $1$ & 2 & 3 & $5$ & 10 & 20 & 48 \\
\hline
$0$           &1.15  & 1.16 & 1.13    & 1.11    & 1.08 & 1.04 & 1.02 & 1.01\\
$0.5$        & 1.02 & 1.03 & 1.05   & 1.06    & 1.08 & 1.10 &1.10 & 1.12 \\
$1$           & 1.004 & 1.02  & 1.06 & 1.09 & 1.14 & 1.20 &1.26 & 1.30 \\
$3$           & 1.03 & 1.06 & 1.13    & 1.20   & 1.33 & 1.55 & 1.79 & 2.03 \\
$10$         & 1.04 &  1.08 &  1.17  & 1.26  & 1.43 & 1.84 & 2.50 & 3.58 \\
\hline
\end{tabular}
\caption{The limit risk ratio $\rho_{\mathcal{E}}(\alpha,\beta)$ for several 
values of  $\alpha, \beta$.}
\label{table:1}
\end{table}

Clearly, $\rho_1(\bar{n}_u,\bar{n}_o,\mathcal{E}) \ge 1$ always holds and 
the bigger $\rho_1$, the bigger the improvement provided by the optimal 
re-allocation.
The value $\rho_1(\bar{n}_u,\bar{n}_s,\mathcal{E})$ gives of course  
a lower bound for the performance of the optimal allocation 
$\bar{n}_o(B_e(\bar{n}_u(n), \mathcal{E})$ defined by Theorem \ref{th1}:
\[
\rho_1(\bar{n}_u,\bar{n}_o,\mathcal{E}) \ge 
\rho_1(\bar{n}_u,\bar{n}_s,\mathcal{E}).
\]
For example, the first row of Table \ref{table:1} gives the relative 
improvement (for several values of $\alpha$) 
of the Pinsker constant provided by the sub-optimal allocation 
$\bar{n}_s$. The optimal allocation $\bar{n}_o$ will do better, not worse at least.

A numerical inspection of the limit risk ratio 
$\rho_{\mathcal{E}}(\alpha,\beta)$ reveals the actual sub-optimality 
of the allocation $\bar{n}_s$. Namely, there is a set $S$ of 
values of $\alpha,\beta$ on which the sub-optimal allocation 
$\bar{n}_s$ fails to improve on $\bar{n}_u$: 
$\rho_\mathcal{E}(\alpha,\beta)<1$ for  $(\alpha,\beta)\in S\subset 
(0,0.3205]\times [0.7,1.823]$. The set $S$ has a ``hill'' form 
and is plotted in Figure \ref{fig1}.
On the positive side, the failure occurs by a relatively small margin: 
$\min_{(\alpha,\beta)\in S}  \rho_\mathcal{E}(\alpha,\beta)={\rho}_o \approx 
\rho_\mathcal{E}(0.149,1.079)=0.998477$.  
\begin{figure}[h]
\begin{center}
\includegraphics[height=2in, width=3.5in] 
{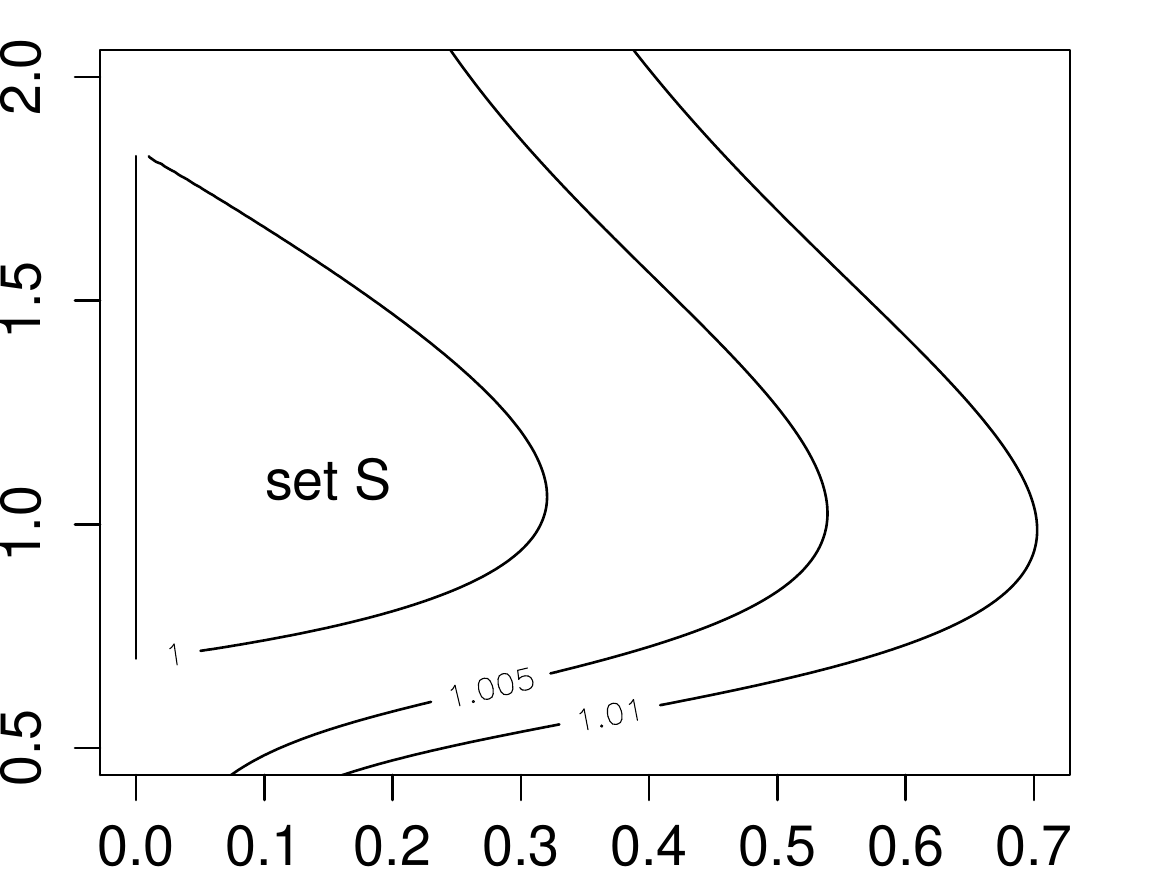}
 \end{center}
 \vspace{-10pt}
\caption{The contour plot of $\rho_{\mathcal{E}}(\alpha,\beta)$; 
the $\alpha$-axis is horizontal.}
\label{fig1}
\end{figure}

\end{example}

\begin{remark}
We conjecture that $\rho_\mathcal{E}(\alpha,\beta)\ge \rho_0\approx 0.998477$ 
for all $\alpha>0$, $\beta>-0.5$. 
Since the beta function is present in the expression for $\rho_\mathcal{E}$, 
this would implicitly  yield an inequality for the beta function: 
for all $\alpha>0$, $\beta>-0.5$,
\[
B^2\Big(\frac{\beta+1}{\alpha},\frac{1}{2}\Big) \le  
\frac{\rho_o^{-1}\alpha(\alpha+2\beta+2)^2}{(2\beta+1)(\alpha+2\beta+1)}  
\bigg(\frac{\rho_o^{-1}2(\beta+1)(2\alpha+2\beta+1)}{(2\beta+1)(3\alpha+2\beta+2)}
\bigg)^{\frac{\beta+1}{\alpha}}.
\]
Notice that the left hand side is a function of $\frac{\beta+1}{\alpha}$, whereas
the right hand side is a function of $\alpha, \beta$ and cannot be 
manipulated into a function of $\frac{\beta+1}{\alpha}$; we say
the above inequality is \emph{unfolded} in $\alpha,\beta$.
\end{remark}

\section{Optimal allocation for hyperrectangles}
\label{sec_hyperrectangles}

Consider $\Theta=\mathcal{H}=\mathcal{H}(a)$ defined in \eqref{ellipsoids} 
and assume for simplicity 
that $\sigma_i$ is nondecreasing with $i$.
First derive the minimax linear risk $R(\bar{n}, \mathcal{H})$ defined 
by \eqref{lin_minimax}. It is easy to see (cf.\  Belitser (2001)) that 
\begin{align}
R(\bar{n},\mathcal{H})&= 
\inf_\lambda\sup_{\theta \in \mathcal{H}} 
\mathrm{E}_\theta \|\hat{\theta}(\lambda) - \theta\|^2
=\sup_{\theta \in \mathcal{H}} 
\inf_\lambda\mathrm{E}_\theta \|\hat{\theta}(\lambda) - \theta\|^2 \notag\\
\label{minimax_hyper}
&=
\sup_{\theta \in \mathcal{H}}  
\mathrm{E}_\theta \|\hat{\theta}(\lambda_o) - \theta\|^2
=\sup_{\theta \in \mathcal{H}}
\sum_i \frac{\theta_i^2 \frac{\sigma_i^2}{n_i}}{\theta_i^2+\frac{\sigma_i^2}{n_i}}
=\sum_i \frac{\sigma_i^2}{n_i+\sigma_i^2a_i^{-2}},
\end{align}
where $\hat{\theta}_i(\lambda_o) = \lambda_{o,i} X_i$ with
$\lambda_{o,i} = \frac{n_i a_i^2}{n_i a_i^2+\sigma_i^2}$.
Minimizing  $R(\bar{n}, \mathcal{H})$ under the restriction $\sum_i n_i =n$
by Lagrange multiplier yields the optimal allocation $\bar{n}_o(n)$:
\begin{align}
\label{n_o_hyper}
n_{o,i}&= \frac{\sigma_i (n+\sum_{k=1}^{d_o} 
\sigma_k^2 a_k^{-2})}
{\sum_{k=1}^{d_o} \sigma_k} -\sigma_i^2a_i^{-2},
\quad i=1,\ldots,d_o,
\end{align}
and $n_{o,i}=0$ for $i >d_o$. Here 
$d_o=d_o(n)=d_o(n,\mathcal{H})$ is the number of nonzero coordinates 
in $\bar{n}_o(n)$ and it is found as the biggest natural 
number such that all $n_{o,i}$ in \eqref{n_o_hyper} are nonnegative. 
In view of monotonicity of $\sigma_i a_i^{-2}$, 
the explicit formula for $d_o$ is readily obtained:
\begin{align}       
\label{d_o_hyper}
d_o=\max\Big\{k \in \mathbb{N}: \sum_{i=1}^k 
\sigma_i(\sigma_ka_k^{-2} -\sigma_i a_i^{-2}) \le n \Big\}.
\end{align}
The optimized risk is then 
\begin{align}
\label{opt_risk_hyper}
R(\bar{n}_o,\mathcal{H}) =
\inf_{\bar{n}\in\mathcal{N}_n} R(\bar{n},\mathcal{H})= 
\frac{\big[\sum_{i=1}^{d_o} \sigma_i \big]^2}
{n+\sum_{i=1}^{d_o} \sigma_i^2a_i^{-2}} +\sum_{i=d_o+1}^\infty a_i^2.
\end{align}
Let us summarize the obtained result by the following theorem.
\begin{theorem}
\label{th2}
The optimal allocation $\bar{n}_o(n)$ is defined by \eqref{n_o_hyper} with 
$d_o$ defined by \eqref{d_o_hyper} and the corresponding optimized 
risk $R(\bar{n}_o,\mathcal{H}(a))$ is defined  by \eqref{opt_risk_hyper}.
\end{theorem}

\begin{remark}
In general case when the sequence $\sigma_i a_i^{-2}$ is not necessarily monotone, 
the formula for $\bar{n}_o(n)$ is a bit more complicated:
\[
n_{o,i}= \frac{\sigma_i (n+\sum_{k\in \mathcal{M}_o} 
\sigma_k^2 a_k^{-2})}
{\sum_{k\in\mathcal{M}_o} \sigma_k} -\sigma_i^2a_i^{-2},
\quad i\in \mathcal{M}_o,
\]
and $n_{o,i}=0$ for $i \not\in \mathcal{M}_o$, where the set 
$\mathcal{M}_o$ is such that $\sum_{k\in\mathcal{M}_o} 
\sigma_k(\sigma_ia_i^{-2} -\sigma_k a_k^{-2}) \le n$ for all $i \in \mathcal{M}_o$.
\end{remark}

\begin{example} 
\label{example2}
Consider the example of Sobolev hyperrectangle 
$\mathcal{H}(a)$ with $a_i^2=Qi^{-(2\alpha+1)}$, 
$\sigma_i^2=\sigma^2 i^{2\beta}$,  $\alpha>0$, $\beta>-\frac{1}{2}$;
$Q=\sigma^2=1$ without loss of generality.
According to \eqref{n_o_hyper}, we derive the optimal allocation $\bar{n}_o(n)$:
\[
n_{o,i}=\frac{n i^{\beta}+
i^\beta \sum_{k=1}^{d_o} k^{2\alpha+2\beta+1}}{\sum_{k=1}^{d_o} k^\beta}
-  i^{2\alpha+2\beta+1}, \quad i=1,\ldots, d_o,
\]
where $d_o=d_o(n)$ is the biggest number $k\in\mathbb{N}$ 
for which, according to \eqref{d_o_hyper},
\[
k^{2\alpha+\beta+1}\sum_{i=1}^k i^{\beta}-
\sum_{i=1}^k i^{2\alpha+2\beta+1}\le n.
\]
Using the relation $\sum_{k=1}^M k^\kappa\sim \frac{M^{k+1}}{\kappa+1}$
as $M \to \infty$ for $\kappa >-1$, we derive  
\[
d_o=d_o(n)\sim B_o n^{\frac{1}{2\alpha+2\beta+2}},  
\;
B_o=B_o(\alpha,\beta) 
=\Big(\frac{2(\beta+1)(\alpha+\beta+1)}{2\alpha+\beta+1}\Big)^{\frac{1}{2\alpha+2\beta+2}}.
\] 
and the optimized risk is
\begin{align*}
R(\bar{n}_o(n),\mathcal{H}) &=
\frac{\big[\sum_{i=1}^{d_o} i^{\beta}\big]^2}
{n+\sum_{i=1}^{d_o}i^{2\alpha+2\beta+1}}
+\sum_{i=d_o+1}^\infty i^{-(2\alpha+1)}
\sim B_{\mathcal{H}}n^{-\frac{\alpha}{\alpha+\beta+1}}.  
\end{align*}
where $B_{\mathcal{H}} =B_{\mathcal{H}}(\alpha,\beta) =
\frac{2\alpha+\beta+1}{2\alpha(\beta+1)} 
\big(\frac{2\alpha+\beta+1}{2(\beta+1)(\alpha+\beta+1)}
\big)^{\frac{\alpha}{\alpha+\beta+1}}$.
In view of \eqref{minimax_hyper}, the minimax linear risk for the uniform 
allocation is  
\[
R(\bar{n}_u(n),\mathcal{H}) = 
\sum_i \frac{i^{2\beta}}{n+i^{2\alpha+2\beta+1}}
\sim\bar{B}_{\mathcal{H}} n^{-\frac{2\alpha}{2\alpha+2\beta+1}},
\]
with (cf.\  Belitser (2001)) $\bar{B}_{\mathcal{H}}=\bar{B}_{\mathcal{H}}(\alpha,\beta)
= \frac{B\big(\frac{2\alpha}{2\alpha+2\beta+1},\frac{2\beta+1}
{2\alpha+2\beta+1}\big)}{2\alpha+2\beta+1}$,
where $B(\cdot,\cdot)$ is the Beta function.
Note that the last two relations imply that
\[
R(\bar{n}_o(n),\mathcal{H})\asymp n^{-\frac{\alpha}{\alpha+\beta+1}} 
\gg n^{-\frac{2\alpha}{2\alpha+2\beta+1}} \asymp R(\bar{n}_u(n),\mathcal{H}),
\] 
for the same reason as for the Sobolev ellipsoids: the measurement budgets 
are very different (in favor of $\bar{n}_u$) as $B(\bar{n}_u)=\infty$, but $B(\bar{n}_o(n))\le n$.

Let us quantify the amount of 
improvement on the uniform allocation $\bar{n}_u$ provided
by the optimal (re-)allocation \eqref{n_o_hyper}. 
The criterion \eqref{risk_ratio} is not suitable for this purpose, since,
in view of \eqref{minimax_hyper}, the effective measurement budget 
of $\bar{n}_u(k)$ is infinite for hyperrectangles: 
$B_e(\bar{n}_u(k),\mathcal{H})=\infty$.  

We use the criterion \eqref{risk_ratio1} instead. The numerator in 
 \eqref{risk_ratio1} is found to be
\[
\inf_{\bar{\delta}, k: \, \bar{\delta}\in \mathcal{B}(\bar{n}_u(k),\bar{n}_o(n))}
R(\bar{n}_u(k)\cdot \bar{\delta},  \mathcal{H})
= \inf_{d,k:\, kd \le n} R(\bar{n}_u(k) \cdot \bar{1}_d,\mathcal{H}),
\vspace{-10pt}
\]
where the \emph{truncation pattern} $\bar 1_d=(\overbrace{1,\ldots, 1}^{d}, 0,0, \ldots)$ 
leads to the truncated version of $\bar{n}_u(k)$:
$\bar{n}_{ut}(k,d)=\bar{n}_u(k)\cdot \bar{1}_d=(k,\ldots, k, 0,0,\ldots)$. 
Using this, by some tedious computations, we derive the numerator in \eqref{risk_ratio1}:  
\begin{align*}
\inf_{d,k:\, kd \le n} R(\bar{n}_u(k) \cdot \bar{1}_d,\mathcal{H})
&=\inf_{d,k:\, k d \le n} R(\bar{n}_{ut}(k,d),\mathcal{H}) 
\sim B'_\mathcal{H} n^{-\frac{\alpha}{\alpha+\beta+1}},
\end{align*}
where $B'_\mathcal{H} =\frac{\alpha+\beta+1}{\alpha} 
\Big(\frac{2\alpha\bar{B}_\mathcal{H}}{2\alpha+2\beta+1}
\Big)^{\frac{2\alpha+2\beta+1}{2\alpha+2\beta+2}}$.
The criterion \eqref{risk_ratio1} becomes
\begin{align*}
&\rho_2(\bar{n}_u,\bar{n}_o,\mathcal{H})
=\frac{\inf_{d,k:\, k d \le n} R(\bar{n}_{ut}(k,d),\mathcal{H})}{R(\bar{n}_o(n),\mathcal{H})} 
\sim \frac{B'_\mathcal{H}}{B_\mathcal{H}} =
\rho_{\mathcal{H}}(\alpha,\beta)  \\
&= \Big(\frac{2(\beta+1)(\alpha+\beta+1)}{2\alpha+\beta+1}\Big)^{\frac{2\alpha+\beta+1}{\alpha+\beta+1}}
\bigg( \frac{2\alpha B\big(\frac{2\alpha}{2\alpha+2\beta+1}, \frac{2\beta+1}{2\alpha+2\beta+1} \big)}
{(2\alpha+2\beta+1)^2}\bigg)^{\frac{2\alpha+2\beta+1}{2(\alpha+\beta+1)}}.
\end{align*}

Table \ref{table:2} presents a small selection of  the computed limit ratio  
$\rho_{\mathcal{H}}(\alpha,\beta)$ for several values of parameters 
$\alpha$ and $\beta$.  
\begin{table}[ht]
\centering
\begin{tabular}{| c | c  c c c c c c c |} 
\hline
\multirow{2}*{$\beta$} 
& \multicolumn{8}{ c|}{$\alpha$}   \\ \cline{2-9}                
 & $0.5$ & $1$ & 2 & 3 & $5$ & 10 & 20 & 48 \\
\hline
$0$    &1.46  & 1.31 & 1.19 & 1.14  & 1.09 & 1.05 & 1.02 & 1.01\\
$0.5$ & 1.51 & 1.38 & 1.27 & 1.22  & 1.18 & 1.15 &1.14 & 1.13  \\
$1$    & 1.52 & 1.43  & 1.35 & 1.33 & 1.31 & 1.31 &1.32 & 1.33 \\
$3$    & 1.44 & 1.45  & 1.49  & 1.54 & 1.63 & 1.79 & 1.96 & 2.12 \\
$10$ & 1.26 &  1.31 &  1.42 & 1.52 & 1.73 & 2.18 & 2.87 & 3.91 \\
\hline
\end{tabular}
\caption{The limit risk ratio $\rho_{\mathcal{H}}(\alpha,\beta)$ for several values of  
$\alpha, \beta$.}
\label{table:2}
\end{table}
\end{example}

\begin{remark}
The relation $\rho_2(\bar{n}_u,\bar{n}_o,\mathcal{H}) \ge 1$ 
implies another (like in Example \ref{example1})
unfolded inequality for the beta function: 
for all $\alpha>0$, $\beta>-0.5$,
\[
B\Big(\frac{2\alpha}{2\alpha+2\beta+1}, \frac{2\beta+1}{2\alpha+2\beta+1}\Big)
\ge 
\frac{1}{2\alpha} \bigg(\frac{(2\alpha+2\beta+1)(2\alpha+\beta+1)}
{2(\beta+1)(\alpha+\beta+1)}\bigg)^2.
\]
\end{remark}

\section{Discussion}
\label{sec_discussion}

We conclude this note with a 
a discussion on possible extensions of the  considered problem. 
The reader is invited to elaborate on these. 
\vspace{-\baselineskip}
\paragraph{Other kind of ill-posedness, ellipsoids, hyperrectangles.} 
One can use Theorems \ref{th1} and \ref{th2} to compute the exact 
asymptotic expressions for other kinds of ill-posedness (sequence $\sigma_i^2$) 
in the model and other kind of ellipsoids $\mathcal{E}(a)$ and 
hyperrectangles $\mathcal{H}(a)$. 
For example, severely ill-posed problem with Sobolev ellipsoids or hyperrectangles: 
$\sigma_i^2=\sigma^2 e^{2\beta i}$ and $a^2_i =Qi^{-(2\alpha+1)}$; 
mildly ill-posed problem with analytic ellipsoids or hyperrectangles:
$\sigma_i^2=\sigma^2 i^{2\beta}$ and $a^2_i =Qe^{-2\alpha i}$.
Other combinations  of ill-posedness and ellipsoids 
or hyperrectangles can be considered.  
\vspace{-\baselineskip} 
\paragraph{Other classes $\Theta$.}
Other choices for the set $\Theta$ can be considered: tail classes,
parametric classes, $\ell_p$-bodies and balls, also sparsity classes.
Consider, for example, one sparsity class: \emph{nearly black vectors}
$\theta\in  \ell_0[p_n] \subset \mathbb{R}^{N_n}$, where
$0\le p_n=o(N_n)$ and, with $\# S$ denoting the number of elements in $S$,
\begin{align*}
\ell_0[p_n]=\big\{\theta\in\mathbb{R}^{N_n}: \# I(\theta) \le p_n\big\}, \quad
I(\theta)= \{i: \theta_i\not=0\} \subseteq \{1,\ldots, N_n\}.
\end{align*}
For the direct problem $\sigma^2_i=1$ with $p_n, N_n \to \infty$ 
as $n\to \infty$, the minimax risk is known to be 
\[
R(\bar{n}_u(n), \ell_0[p_n]) =\inf_{\hat{\theta}}\sup_{\theta \in\ell_0[p_n]} 
\mathrm{E}_\theta \|\hat{\theta}-\theta\|^2
\sim \frac{2p_n \log N_n}{n}.
\]
As to the indirect case $\sigma_1^2 \ge \sigma_2^2 \ge \ldots \ge \sigma_{N_n}^2 >0$ and 
an arbitrary allocation $\bar{n}$, we are unaware of a 
result on the minimax risk, but we conjecture that  
\[
R(\bar{n}, \ell_0[p_n])  
\sim 2 \log N_n \sum_{i \le p_n} \frac{\sigma_i^2}{n_i} \triangleq R_c(\bar{n}, \ell_0[p_n]).
\]
Minimizing the conjectured asymptotic risk $R_c(\bar{n}, \ell_0[p_n])$  
with respect to $\bar{n} \in \mathcal{N}_n$,
we obtain that $\inf_{\bar{n} \in \mathcal{N}_n} R_c(\bar{n}, \ell_0[p_n]) = 
R_c(\bar{n}_o(n), \ell_0[p_n])$ with 
\[
R_c(\bar{n}, \ell_0[p_n])
\sim 2 \log N_n \frac{\big(\sum_{i \le p_n} \sigma_i\big)^2}{n}, \quad 
n_{o,i} = \frac{ n \sigma_i 1\{i\le p_n\}}{\sum_{i \le p_n} \sigma_i}.
\]
If $\sigma_i^2 =1$, $R_c(\bar{n}, \ell_0[p_n]) = \frac{2p_n^2 \log N_n}{n}$ and 
it is easy to see that the optimal allocation does not improve on the 
uniform allocation in the direct case.
\vspace{-\baselineskip}
\paragraph{Adaptive optimal allocation.}
Another interesting extension would be to obtain the adaptive 
versions of the results on the optimal allocation problem, i.e., without knowledge of 
the structural parameter of the set $\Theta$. For example, construct the 
optimal allocation for the Sobolev ellipsoid $\mathcal{E}(a)$ 
(or hyperrectangle), $a_i^2=Qi^{-2\alpha}$, without using 
the smoothness parameter $\alpha$.

\end{document}